\documentclass[11pt]{article}

\usepackage{fancyhdr}
\usepackage{ae}
\usepackage{aecompl}
\usepackage[T1]{fontenc}
\usepackage{graphics}
\usepackage[dvips]{graphicx}
\usepackage{psfrag}
\usepackage[numbers]{natbib}
\usepackage{float}
\usepackage[centertags]{amsmath}
\usepackage{amsfonts}
\usepackage{amssymb}
\usepackage{subfigure}
\usepackage{tabularx}
\usepackage{theorem}
\usepackage{pifont}

\newtheorem{theorem}{Theorem}

\newtheorem{remark}{Remark}
\newtheorem{proposition}{Proposition}
\newtheorem{lemma}{Lemma}

   \title {\bf Analysis of a class of non linear subdivision schemes and associated   multi-resolution
transforms\footnote{Research partially
   supported by European network ``Breaking complexity''
  \# HPRN-CT-2002-00286}
  }

\author{S. Amat\thanks{
 Departamento de Matem\'atica Aplicada y Estad\'istica.  Universidad
 Polit\'ecnica  de Cartagena  (Spain).Research supported in part by  00675/PI/04 and MTM2004-07114.
e-mail:{\tt sergio.amat@upct.es}}\and K. Dadourian\thanks{ Laboratoire
d'Analyse Topologie et Probabilites (LATP), Universit\'e de
Provence (France).
 e-mail:{\tt
dadouria@cmi.univ-mrs.fr}}\and J. Liandrat\thanks{ Ecole
Centrale de  Marseille (ECM). Laboratoire
d'Analyse Topologie et Probabilites (LATP).
 e-mail:{\tt
jacques.liandrat@ec-marseille.fr}} }
\begin{document}

\maketitle
\underline{Corresponding author:} Jacques Liandrat

\vskip1cm


    \begin{abstract}
      This paper is devoted to the convergence and stability analysis of
a class of nonlinear subdivision schemes and associated multi-resolution transforms.
These schemes are defined  as a perturbation of a
linear subdivision scheme. Assuming  a contractivity property,  stability and
convergence  are derived. These results are then applied to various schemes such as
uncentered interpolatory  linear scheme, WENO scheme \cite{LOC}, Power-P scheme \cite{SM} and a non linear scheme using local spherical coordinates \cite{AEV}.
    \end{abstract}


\vspace{10pt}
{\bf AMS(MOS) subject classifications.} 41A05, 41A10, 65D05, 65D17

    \section{Introduction}

\indent Multi-resolution representations of discrete data are useful tools in several
areas of application as image compression or adaptive methods for
partial differential equations. In these applications, the ability of these
 representations to approximate the input data with  high accuracy using
a very small set of coefficients is a central property. Moreover, the stability of these representations in presence of perturbations (generated by compression or due to approximations) is a key point.\\
\indent In the last decade, several attempts to improve the property of  classical linear
multi-resolutions have lead to nonlinear multi-resolutions. In many cases, this
nonlinear nature hinders the proofs of convergence and stability.\\
\indent In \cite{ADLT}, in the context of image compression, a new
multi-resolution transform  has been presented. This multi-resolution is based
on an univariate nonlinear multi-resolution called PPH multi-resolution
(see \cite{KD} in the context of convexity preserving). It has
been a\-na\-ly\-zed in terms of convergence and
stability of an associated subdivision scheme following an
approach for data dependent multi-resolutions introduced in
\cite{CDM}. Due to nonlinearity, the stability of the PPH multi-resolution is not a consequence of the convergence of the associated subdivision scheme. It has been established in \cite{AL},
presenting  the PPH subdivision scheme  as some
perturbation of a  linear scheme following
\cite{DY}, \cite{Os}, \cite{DRS} and \cite{FM}.\\
\indent The aim of the present paper is to generalize the results presented
in \cite{AL} for a general family of nonlinear multi-resolution
schemes associated to an interpolatory subdivision scheme  $S_{NL}: l^{\infty}(\mathbb{R}) \rightarrow l^{\infty}(\mathbb{R})$ of the form:
\begin{equation}
\forall f \in l^{\infty}(\mathbb{R}), \quad \forall n\in \mathbb{Z} \qquad
\left \lbrace \begin{array}{lll}
S_{NL}(f)_{2n+1} $=$ S(f)_{2n+1}+F(\delta f)_{2n+1},\\
S_{NL}(f)_{2n} $=$ f_{n},
\end{array} \right .
\end{equation}
where $F$ is a nonlinear operator defined on $l^{\infty}({\mathbb{R}})$,
$\delta$ is a  linear  and continuous operator
 on $l^{\infty}(\mathbb{R})$ and $S$ is a linear and convergent subdivision scheme. Considering two subdivision schemes $S_{NL}$ and $S$, it is always possible to introduce the difference $F=S_{NL}-S$. If one assume some properties of polynomial reproduction (see section \ref{sec3}), as shown in \cite{KD}, $F$ is in fact  a function of differences, i.e. of $df$ defined by $df_n=f_{n+1}-f_n$.

 \vskip1cm

Theorems \ref{th1} and  \ref{th5}, that are the main results of this paper, establish that if $F,S$
and $\delta$ satisfy some natural properties, then the subdivision scheme is
convergent and the multi-resolution is stable.\\
\indent The paper is organized as follows.
 In section \ref{sec2} we recall briefly the Harten's interpolatory multi-resolution
  framework which is the natural setting  for our work.
 We  precise  the class of  schemes under consideration and we establish the
 main results  in section \ref{sec3}.
Various applications are presented  in section \ref{sec4}.


\section{Harten's framework and basic definitions}\label{sec2}

\indent   In Harten's interpolatory multi-resolution, one considers a set of nested bi-infinite
 regular grids:
  $$X^{j}=\{ x_{n}^{j} \}_{n\in \mathbb{Z}},\quad x_{n}^{j}=n2^{-j},
  $$
  where  $j$ stands for a scale parameter and $n$ controls the position. \\

\indent The point-value discretization operators (or sampling operators) are defined by

\begin{equation} {\cal D}_{j}:  f \in C(\mathbb{R}) \mapsto
f^j=(f_{n}^{j})_{n\in \mathbb{Z}}:=(f(x_{n}^{j}))_{n\in \mathbb{Z}} \in
V^j,
\end{equation}
  where $V^{j}$ is the space of real sequences and $C(\mathbb{R})$ the set of continuous
  functions on $\mathbb{R}$.\\

\indent The  reconstruction operators ${\cal R}_{j}$ associated to  this discretization  are
any right inverses of ${\cal D}_j$ on $V^j$, that is, any operators ${\cal R}_j$ satisfying\,:

\begin{equation}
({\cal R}_{j}f^{j})(x_{n}^{j})=f_{n}^{j}=f(x_{n}^{j}).
\end{equation}

\indent  For any $j$, the operator defined by ${\cal D}_{j}{\cal R}_{j+1}$ acts between
a fine  scale $j+1$ and a coarser scale j. Here, it is  a sub-sampling
operator from $V_{j+1}$ to $V_j$. \\
\indent The operator defined by ${\cal D}_{j+1}{\cal R}_j$ acts between
a coarse scale $j$ and a finer scale $j+1$ and is called a
prediction operator. A prediction operator can be considered as  a subdivision scheme
\cite{Dyn} from $V_j$ to $V_{j+1}$. We say that the
subdivision scheme $S$ defined by  $(f^j) \mapsto S(f^j)={\cal D}_{j+1}{\cal R}_j(f^j)$
 is uniformly convergent if\,:
\begin{equation*}
\forall f \in l^\infty, \exists f^{\infty} \in C^0(\mathbb{R}) \quad
\textrm{ such that} \quad  \displaystyle{\lim_{j\rightarrow +\infty} \sup_{n\in \mathbb{Z}}
}|S^j(f)_n-f^{\infty}(2^{-j}n)|=0.
\end{equation*}
We note $f^\infty=S^\infty f$.

\indent Since for most function $f$, ${\cal D}_{j+1}{\cal R}_j f^j\neq f^{j+1}$, details,
 called $d^{j}$ and defined  by  $d^{j}=f^{j+1}-{\cal D}_{j+1}{\cal R}_j f^j$, should be added to ${\cal D}_{j+1}{\cal R}_j f^j$ to recover $f^{j+1}$ from $f^j$.
The multi-re\-so\-lution decomposition  (see \cite{AD00}, \cite{Har93}, \cite{Har96} for precisions)
 of $f^L$ is the sequence $\{f^0,d^0,\ldots,d^{L-1}  \}$. Moreover,
the multi-resolution transform is said to be stable if\,:
\begin{eqnarray}
\exists   C \mbox{ such that } &  \forall f^L, \tilde{f}^L, j \leq L & \nonumber \\
 & || f^j-\tilde{f}^j||_\infty & \leq C \left (||
f^0-\tilde{f}^0||_\infty + \sum_{k=0}^{j-1} || d^k-\tilde{d}^k||_\infty \right ),
\label{s1}\\
   & || f^0-\tilde{f}^0 ||_\infty & \leq  C || f^j-\tilde{f}^j||_\infty,
\label{s2} \\
  & || d^k-\tilde{d}^k ||_\infty & \leq C  || f^j-\tilde{f}^j ||_\infty,  \quad  \forall  k,\, 0\leq k \leq j-1,
 \label{s3}
\end{eqnarray}
where $\{\tilde{f}^0,\tilde{d}^0,\ldots,\tilde{d}^{L-1}  \}$ is the multi-resolution  decomposition of $\tilde{f}^L$.

 \indent When the prediction operator ${\cal D}_{j+1}{\cal R}_j f^j$ is linear, the
 convergence of the associated subdivision scheme implies the stability of the multi-resolution analysis. In the
 non linear case, it is not the case and there is no general result for the multi-resolution analysis stability.
 \vspace{0.5cm}


\section{A Class of Nonlinear Subdivision Schemes}\label{sec3}

\indent Introducing $S$ a linear, reproducing polynomials
\footnote{  The interpolatory subdivision scheme $S$  reproduces polynomials of
 degree $P$ if, for any polynomials $\cal P$ of degree less or equal to $P$, if
 $f_n={\cal P}(x^j_n)$ then $S(f)_{2n+1}={\cal P}(x^{j+1}_{2n+1})$}  up to degree $P$ and convergent interpolatory subdivision scheme we consider   nonlinear interpolatory subdivision schemes  that write
\begin{equation}
 \left \lbrace \begin{array}{lll}
S_{NL}(f^j)_{2n+1} $=$ S(f^j)_{2n+1}+F(\delta f^j)_{2n+1}\\
S_{NL}(f^j)_{2n} $=$f^j_{n}
\end{array} \right .
\label{75}
\end{equation}
where $F$ is a nonlinear operator defined on $l^{\infty}(\mathbb{Z})$ and  $\delta$ is a continuous
linear operator  on $l^{\infty}(\mathbb{Z})$.
\ \\

\subsection{Convergence analysis}
\label{subconv}

 We have  the following theorem related to the convergence of the nonlinear subdivision scheme $S_{NL}$:\,:

\begin{theorem}

\label{th1} If $F, S $ and $\delta$ verify:

\begin{eqnarray}
 & & \exists M>0 \quad  \textrm{such that} \quad  \forall d \in l^{\infty}
\quad \; || F(d)||_\infty \leq M ||d||_{\infty}\, \label{h1}\\
 & & \exists c<1 \; \textrm{such that} \qquad
|| \delta S(f)+ \delta F(\delta f)  ||_\infty \leq c || \delta f||_\infty,
\label{h2}
\end{eqnarray}
then the subdivision  scheme $S_{NL}$ is uniformly convergent. Moreover, if $S$
is $C^{\alpha^-}$ convergent (i.e for all $f \in l^{\infty}(\mathbb{Z}), S^\infty(f) \in C^{\alpha^-}$ \footnote{For $ 0<\alpha\leq 1$,
 $C^{\alpha^-}=\{f\textrm{ continuous,  bounded and verifying }\forall
  \alpha^{'}<\alpha, \,\exists C>0, \, \forall x,y\in \mathbb{R} , \quad |f(x)- f(y)|\leq C
  |x-y|^{\alpha^{'}} \}$.
  For $\alpha>1$ with $\alpha=p+r>0$, $p \in I\!\!N$ and $0<r<1$,
  $C^{\alpha^-}(\mathbb{R})=\{f \textrm{ with } f^{(p)}\in C^{r^-}\}$
  }) then, for all sequence $f \in l^{\infty} (\mathbb{Z}),\;S^{\infty}_{NL}(f) \in C^{\beta^-}$ with $\beta=\min{\left (\alpha,log_2(c) \right )}$. \\

\end{theorem}

{\bf Proof}
Using hypotheses (\ref{h1}) and (\ref{h2}) and the definition of $S_{NL}$, we get\,:
\begin{eqnarray*}
|S_{NL}(f^j)_{2n+1}-S(f^j)_{2n+1}| & \leq & M \Vert \delta f^j\Vert_\infty,\\
\Vert S_{NL}(f^j)-S(f^j)\Vert_\infty & \leq & M \Vert \delta (S_{NL}f^{j-1})\Vert_\infty, \\
\Vert S_{NL}(f^j)-S(f^j)\Vert_\infty& \leq & Mc\Vert
\delta f^{j-1}\Vert_\infty,
\end{eqnarray*}
that can be rewritten as\,:
\begin{eqnarray*}
\label{eq38} \Vert S_{NL}(f^j)-S(f^j)\Vert_\infty &\leq  &  M
c^{j} \Vert \delta f^{0} \Vert_{\infty}.
\end{eqnarray*}
Writing\,:
\begin{eqnarray}
\label{majtheorem1}
\Vert S_{NL}(f^j)-S(f^j)\Vert_\infty &\leq  &  M \Vert \delta f^{0} \Vert_{\infty}
2^{jlog_2(c)}
\end{eqnarray}
the convergence  of the
subdivision scheme $S_{NL}$ can be obtained  applying theorem 3.3 of \cite{DRS}.\\
In our context, this theorem applies as follows\,: \\


 If $S$ is a linear $C^{\alpha^-}$ convergent subdivision scheme reproducing polynomials up to degree $P$ and if $S_{NL}$ is a perturbation of $S$ in the sense that, calling $f^k:=S_{NL}(f^0)$ for all $f^0 \in l_\infty$,
 $$ ||S_{NL}(f^k)-S(f^k)||_\infty=O(2^{-\nu k}),$$
 then $S_{NL}$ is $C^{\beta^-}$convergent with  $\beta= \geq
 \mbox{min}(P, s_L, \nu)$.

It follows that if  $S$ is $C^{\alpha^-}$ convergent then
$S_{NL}$ is
 at least $C^{\beta^-}$ convergent  with $\beta=\min{\left (\alpha,log_2(c) \right )}$.



\ \\
$\Box$
\ \\
\begin{remark}
When $F$ is linear, theorem \ref{th1} is a consequence of   theorem 6.2 in \cite{Dyn}.
\end{remark}

\begin{remark}
\label{regularite}
In many  of our examples,  $S$ is the two point centered linear   scheme  defined by $S(f^j)_{2n+1}=\frac{f^j_n+f^j_{n+1}}{2}$ which
limit function is in $C^{1-}$.
Therefore, as soon as   the non linear scheme $S_{NL}$  verifies
hypothesis (\ref{h1}) and (\ref{h2}) with $c \geq  \frac{1}{2}$, the $S_{NL}$ is    $C^{(-log_2(c))^-}$ convergent.
\end{remark}

\begin{remark}
 Hypothesis \ref{h2} can be weakened as:
\begin{eqnarray*}
\exists p \in \mathbb{N} &\exists c<1 &  \textrm{such that} \qquad
\Vert \delta (S_{NL}^p f)  \Vert_\infty \leq c \Vert \delta f\Vert_\infty.
\end{eqnarray*}
The proof remains   the same except that\,:
\begin{eqnarray*}
\Vert S_{NL}(f^j)-S(f^j)\Vert_\infty &\leq  &  M
c \Vert \delta f^j \Vert_{\infty}
\end{eqnarray*}
becomes\,:
\begin{eqnarray*}
\Vert S_{NL}(f^j)-S(f^j)\Vert_\infty &\leq  &  M  \Vert\delta (S_{NL}^p f^{j-p})
\Vert_{\infty},
\\
\Vert S_{NL}(f^j)-S(f^j)\Vert_\infty &\leq  &  M c \Vert \delta f^{j-p} \Vert_{\infty},
\end{eqnarray*}
that can be rewritten, for $j\equiv i [p]$, as:
\begin{eqnarray*}
\Vert S_{NL}(f^j)-S(f^j)\Vert_\infty & \leq  &  M
c^{\frac{j-i}{p}} \Vert \delta f^{i} \Vert_{\infty}.
\end{eqnarray*}
The conclusion is reached applying  theorem 3.3 of \cite{DRS}.
\end{remark}
\begin{remark}

\label{deuxoperateurs}
 A straightforward generalization of theorem  \ref{th1} can be obtained introducing two linear operator
$\delta_1$, $\delta_2$ and a perturbation of the form $F(\delta_1 f,\delta_2f)$.
 Under the following hypotheses:
\begin{eqnarray}
\label{deux1}
\exists M >0 \textrm{ such that  }  |F(d,d')|&\leq& M
\max{\left ( ||d||_{\infty},||d'||_{\infty} \right )},\\
\label{deux2}
\exists c>1 \textrm{ such that  }
||\delta_1(S_{NL}(f))||_{\infty}&\leq&c
\max{\left ( ||\delta_1f||_{\infty},||\delta_2f||_{\infty} \right ) },\\
\label{deux3}
  ||\delta_2(S_{NL}(f))||_{\infty}&\leq&c
\max{\left ( ||\delta_1f||_{\infty},||\delta_2f||_{\infty} \right ) },
\end{eqnarray}

for all $d,d'\in l^{\infty},  f\in l^{\infty}$,
the scheme $S_{NL}$ is uniformly convergent.
\end{remark}

\begin{remark}
\label{R2}
We can also apply theorem \ref{th1} to bi-variate  schemes written as
\begin{eqnarray*}
S_{NL}(x^j,y^j) =
\left (
\begin{array}{cc}
S_{NL_1}(x^j,y^j)_{2n+1}\\
S_{NL_2}(x^j,y^j)_{2n+1}
\end{array}\right ).
=
\left (
\begin{array}{cc}
x_{2n+1}^{j+1}\\
y_{2n+1}^{j+1}
\end{array}\right )
= \left (
\begin{array}{cc}
 S(x^j)_{2n+1}+F_1(\delta x^j,\delta y^j)\\
 S(y^j)_{2n+1}+F_1(\delta x^j,\delta y^j)
\end{array}\right )
\end{eqnarray*}
If the following conditions are satisfied for $i=1,2$
\begin{eqnarray*}
\exists M >0 \textrm{ such that  }  |F_i(d,d')|&\leq& M
\max{\left ( ||d||_{\infty},||d'||_{\infty} \right )},\\
\exists c>1 \textrm{ such that  }
||\delta(S_{NL_i}(x,y))||_{\infty}&\leq &c
\max{\left ( ||\delta x||_{\infty},||\delta y||_{\infty} \right ) },
\end{eqnarray*}
for all $d,d',x,y\in l^{\infty}$,
the scheme $S_{NL}$ is uniformly convergent.
\end{remark}

\ \\




%
\subsection{Stability analysis}
\label{substab}

 We now  consider the multi-resolution analysis associated to the subdivision scheme
(\ref{75}) recalling that,  for any sequence $f^j$, the details  $d^j$
are defined  by $d_n^j=f_{2n+1}^{j+1}-S_{NL}(f^{j})_{2n+1}$.\\

\indent We have   the following theorem concerning the stability of the multi-resolution\,:\\

\begin{theorem}\label{th5}

If $F, S$ and $\delta$ verify:
$\exists M>0, c<1$ such that $\forall f,g, d_1,d_2$,
\begin{eqnarray}
  \quad  ||F(d_1)-F(d_2)||_\infty  \leq M ||d_1-d_2||_{\infty}, \label{eq30}\\
  \quad  \Vert \delta (S_{NL}f-S_{NL}g)  \Vert_\infty \leq c \Vert \delta (f-g)
\Vert_\infty, \label{eq31}
\end{eqnarray}
then  the multi-resolution transform associated
to the  non linear subdivision scheme $S_{NL}$ is stable.
\end{theorem}
\ \\
{\bf Proof}\\
\ \\
 We first prove (\ref{s1})\,:\\
Due to the interpolatory property, we only consider
 $|f^j_{2n+1}-\tilde{f}^j_{2n+1}|$.

Since  $S$ is a convergent linear scheme, we have,  using the stability of the linear scheme $S$: $\exists C'>0$ such that
\begin{eqnarray}
\nonumber
|f^j_{2n+1}-\tilde{f}^j_{2n+1}|& \leq &
C' \left (
||f^0-\tilde{f}^0||_\infty+\sum_{k=1}^{j}||f^k-S(f^{k-1})-\tilde{f}^k+S(\tilde{f}^{k-1})||_{\infty} \right ) \\
\nonumber
& \leq &
C' \left ( ||f^0-\tilde{f}^0||_\infty+\sum_{k=1}^{j}||d^{k-1}+F(\delta
f^{k-1})-\tilde{d}^{k-1}-F(\delta
\tilde{f}^{k-1})||_{\infty} \right ). \\
\nonumber
\end{eqnarray}
 From (\ref{eq30})\,:
\begin{eqnarray*}
\nonumber
|f^j_{2n+1}-\tilde{f}^j_{2n+1}|& \leq &
C'\left ( ||f^0-\tilde{f}^0||_\infty+\sum_{k=0}^{j-1}||d^k-\tilde{d}^k||_{\infty}+M
\sum_{k=1}^{j}||\delta(f^{k-1})-\delta(\tilde{f}^{k-1}) ||_{\infty} \right ) . \\
\end{eqnarray*}
Concentrating on  the last right hand side term we get:

\begin{eqnarray*}
 \nonumber
\sum_{k=1}^{j} ||\delta(f^{k-1})-\delta(\tilde{f}^{k-1}) ||_{\infty}& \leq &
\Vert \delta(f^0)-\delta(\tilde{f}^{0}) \Vert_{\infty}\\
 & & +\sum_{k=2}^{j} \left ( \Vert \delta(S_{NL}f^{k-2})-\delta (S_{NL}\tilde{f}^{k-2}) \Vert_{\infty}+ \Vert \delta
d^{k-2}- \delta \tilde{d}^{k-2} \Vert_{\infty} \right )
\end{eqnarray*}
From (\ref{eq31})\ we get:
\begin{eqnarray*}
\nonumber
\sum_{k=1}^{j}||\delta(f^{k-1})-\delta(\tilde{f}^{k-1}) ||_{\infty}& \leq &
\Vert \delta(f^0)-\delta(\tilde{f}^{0}) \Vert_{\infty}
 +\sum_{k=0}^{j-2} \left ( c \Vert \delta(f^{k})-\delta (\tilde{f}^{k}) \Vert_{\infty}+ \Vert \delta
d^{k}- \delta \tilde{d}^{k} \Vert_{\infty} \right )\\
&\leq  & \sum_{k=0}^{j-2} \left (  c^{k} \Vert \delta f^0 -\delta \tilde{f}^0
\Vert_\infty + \sum_{l=0}^{k} c^{k-l} \Vert
\delta d^{l}-\delta \tilde{d}^{l}\Vert_\infty \right ).
\end{eqnarray*}
Since $0<c<1$ we get finally:
\begin{eqnarray*}
\Vert f^{j}-\tilde{f}^j \Vert_\infty & \leq &   C'||f^0-\tilde{f}^0||_\infty
+C'\sum_{k=0}^{j-1}||d^k-\tilde{d}^k||_{\infty}
 \\
 &&+MC'\frac{1}{1-c} \left (||\delta(f^0)-\delta(\tilde{f}^0)||_{\infty}+
 \sum_{k=0}^{j-2} \Vert \delta d^k-\delta
 \tilde{d}^k\Vert_\infty \right ), \\
\end{eqnarray*}
and, using the continuity of  $\delta $, we get (\ref{s1}) with
a constant
\begin{equation*}
C=C'+\frac{MC' \Vert \delta \Vert_{\infty}}{1-c}.
\end{equation*}
\ \\
We now establish (\ref{s2}) et (\ref{s3}).\\
Equation  (\ref{s2}) is a direct consequence of the interpolatory properties.\\
For (\ref{s3}),  we have,  for $0\leq k \leq j-1$\,:
\begin{eqnarray*}
\nonumber
|d_{n}^k-\tilde{d}_{n}^k| &\leq & ||f^{k+1}-\tilde{f}^{k+1}-S(f^{k})-S(\tilde{f}^{k})||_{\infty} +
||F(\delta
f^{k})-F(\delta \tilde{f}^{k})||_{\infty}.
\end{eqnarray*}
Using the property (\ref{s3}) for the multi-resolution associated to $S$, hypothesis (\ref{eq30}) and the continuity of $\delta$, we have\,:
\begin{eqnarray*}
|d_{n}^k-\tilde{d}_{n}^k| &\leq & C' || f^{j}-\tilde{f}^{j}||_\infty + M \Vert \delta
\Vert_\infty  \Vert  f^{k}- \tilde{f}^{k} \Vert_\infty. \\
\end{eqnarray*}
From (\ref{s2}) for the multi-resolution associated to $S_{NL}$\, we have
\begin{eqnarray*}
|d_{n}^k-\tilde{d}_{n}^k| &\leq &  C' || f^{j}-\tilde{f}^{j}||_\infty + M \Vert \delta
\Vert_\infty  \Vert  f^{j-1}- \tilde{f}^{j-1} \Vert_\infty, \\
\end{eqnarray*}
and  therefore we get (\ref{s3}) with $C=C' + M||\delta||_\infty $. \\
$\Box$\\
\ \\
\begin{remark}

As previously, we can again consider a weaker formulation for hypothesis (\ref{eq31}) such as:
\begin{eqnarray*}
 \exists p \in \mathbb{N}, & \exists c<1 & \textrm{such that} \quad  \Vert \delta (S_{NL}^p f-S_{NL}^p g)  \Vert_\infty \leq c \Vert \delta (f-g)
\Vert_\infty .
\end{eqnarray*}
Under this hypothesis, the stability of the subdivision scheme can still be established. However,  the multi-resolution stability is not ensured. To get it, a stronger hypothesis like:

\begin{eqnarray*}
\exists p \in \mathbb{N}, \exists c<1, \textrm{such that} & &  \\
\quad  \Vert \delta (f^p f-g^p )  \Vert_\infty \leq c \Vert \delta (f-g)
\Vert_\infty&+&M \sum_{k=0}^{p-2} ||d^k(f)-d^k(g)||_{\infty},
\end{eqnarray*}

is required.
\end{remark}

\ \\
\section{Applications}\label{sec4}

\indent This sections is devoted to applications of the previous results to three specific
subdivision schemes (linear and nonlinear) available in the literature. We  provide for each of them, the proofs of convergence and stability.

In all this section, given $f^j=(f^j_k)_{k \in Z\!\!Z}$ we note:

\begin{eqnarray}
df^j&=&(df^j_n)_{n \in Z\!\!Z} \mbox{ with } df^j_n = f^j_{n+1}-f^j_n,\\
Df^j&=&(Df^j_n)_{n \in Z\!\!Z}  \mbox{ with } Df^j_n =  f^j_{n+1}-2f^j_n+f^j_{n-1},\\
\end{eqnarray}

and, more generally $D^lf^j=(D^lf^j_n)_{n \in Z\!\!Z}$ with:

\begin{eqnarray}
D^lf^j_n = D(D^{l-1}f^j)_n=\sum_{i=0}^{2l}(-1)^{i}C_{2l}^ if_{n-l+i} &\textrm{ with }&
 C^i_k=\frac{k!}{i!(k-i)!}.
\end{eqnarray}

\subsection{Multi-resolution analysis associated to a linear  fully non centered Lagrange interpolatory subdivision scheme}

\indent As it has been said before,   for linear scheme, the stability of the multi-resolution analysis is a consequence of
 the convergence of the subdivision scheme (see \cite{Har96} ). Therefore, we only  consider here the convergence of the subdivision scheme. \\
\indent The convergence of centered linear interpolatory schemes is well
know since Delauriers and Debuc \cite{DeDu}.\\ For linear but non centered schemes there is no general results of convergence. Moreover the general tools proposed in \cite{Dyn} are very fastidious to apply and don't provide general results.

\begin{figure}[h!]
\center
\subfigure[9 point scheme]
{\psfrag{0}[l][l]{\footnotesize{$0$ }}
 \psfrag{0.5}[l][l]{\footnotesize{$0.5$ }}
 \psfrag{1}[l][l]{\footnotesize{$1$ }}
 \psfrag{-3}[l][l]{\footnotesize{$-3$ }}
 \psfrag{-2}[l][l]{\footnotesize{$-2$ }}
 \psfrag{-1}[l][l]{\footnotesize{$-1$ }}
\psfrag{0}[l][l]{\footnotesize{$0$ }}
\psfrag{1}[l][l]{\footnotesize{$1$ }}
\psfrag{2}[l][l]{\footnotesize{$2$ }}
\psfrag{3}[l][l]{\footnotesize{$3$ }}
\includegraphics[width=0.45\textwidth]{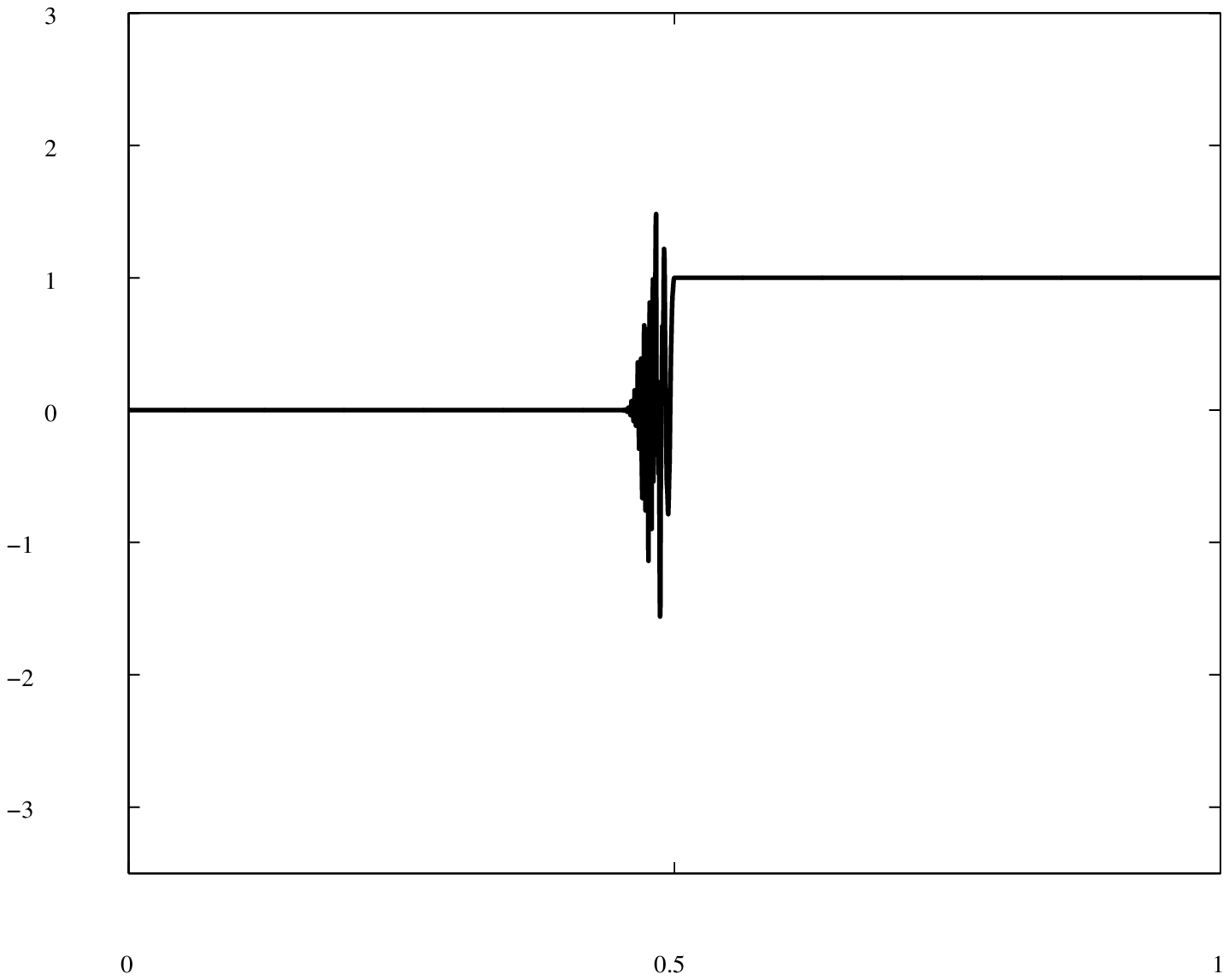}
}
\subfigure[10 point scheme ]
{ \psfrag{0}[l][l]{\footnotesize{$0$ }}
 \psfrag{0.5}[l][l]{\footnotesize{$0.5$ }}
 \psfrag{1}[l][l]{\footnotesize{$1$ }}
 \psfrag{-3}[l][l]{\footnotesize{$-3$ }}
 \psfrag{-2}[l][l]{\footnotesize{$-2$ }}
 \psfrag{-1}[l][l]{\footnotesize{$-1$ }}
\psfrag{0}[l][l]{\footnotesize{$0$ }}
\psfrag{1}[l][l]{\footnotesize{$1$ }}
\psfrag{2}[l][l]{\footnotesize{$2$ }}
\psfrag{3}[l][l]{\footnotesize{$3$ }}
\includegraphics[width=0.45\textwidth]{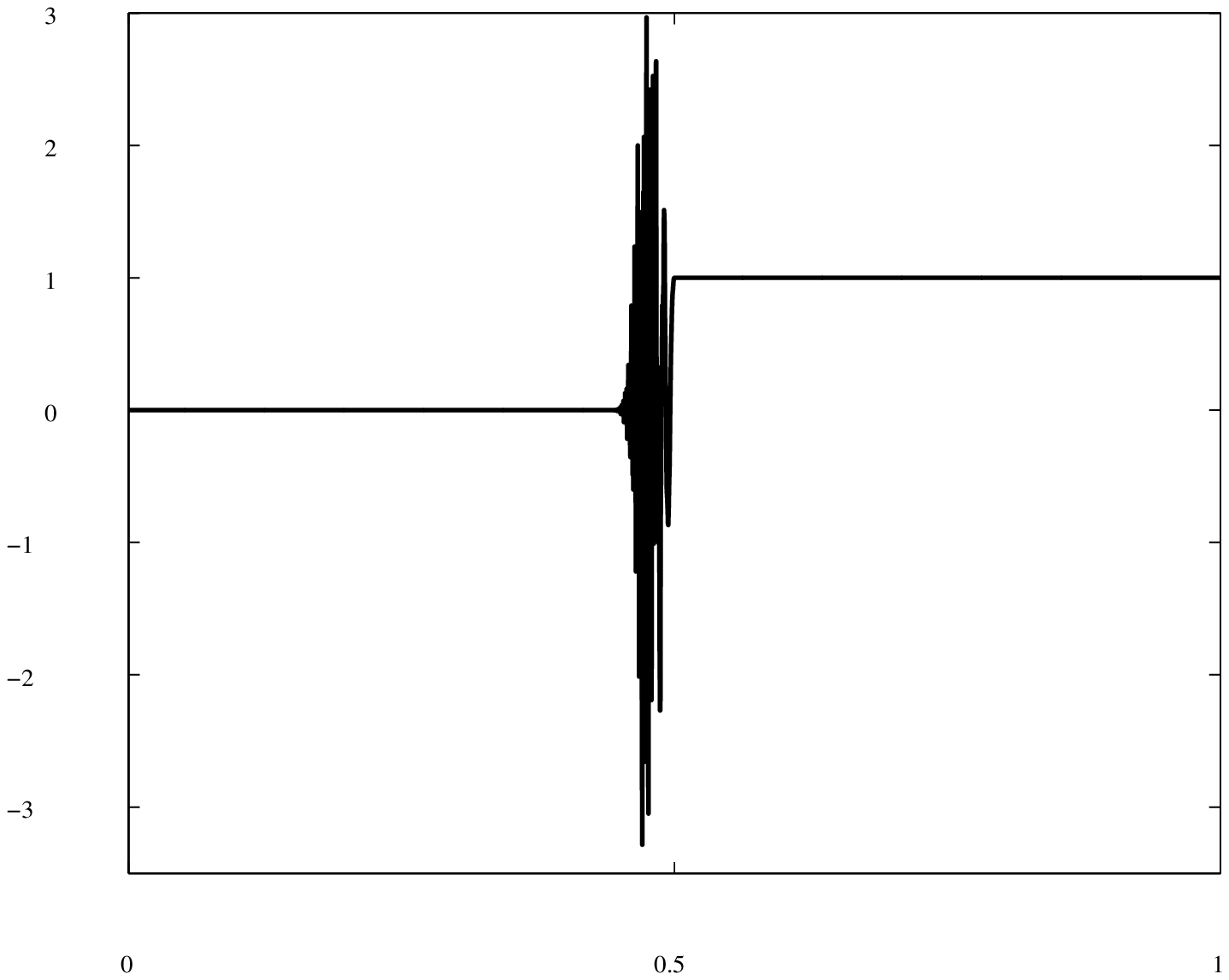}
 }
\\
 \subfigure[9 point scheme (zooming)]
{\psfrag{0.45}[l][l]{\footnotesize{$0.45$ }}
 \psfrag{0.5}[l][l]{\footnotesize{$0.5$ }}
 \psfrag{0.47}[l][l]{\footnotesize{$0.47$ }}
 \psfrag{-3}[l][l]{\footnotesize{$-3$ }}
 \psfrag{-2}[l][l]{\footnotesize{$-2$ }}
 \psfrag{-1}[l][l]{\footnotesize{$-1$ }}
\psfrag{0}[l][l]{\footnotesize{$0$ }}
\psfrag{1}[l][l]{\footnotesize{$1$ }}
\psfrag{2}[l][l]{\footnotesize{$2$ }}
\psfrag{3}[l][l]{\footnotesize{$3$ }}
 \includegraphics[width=0.45\textwidth]{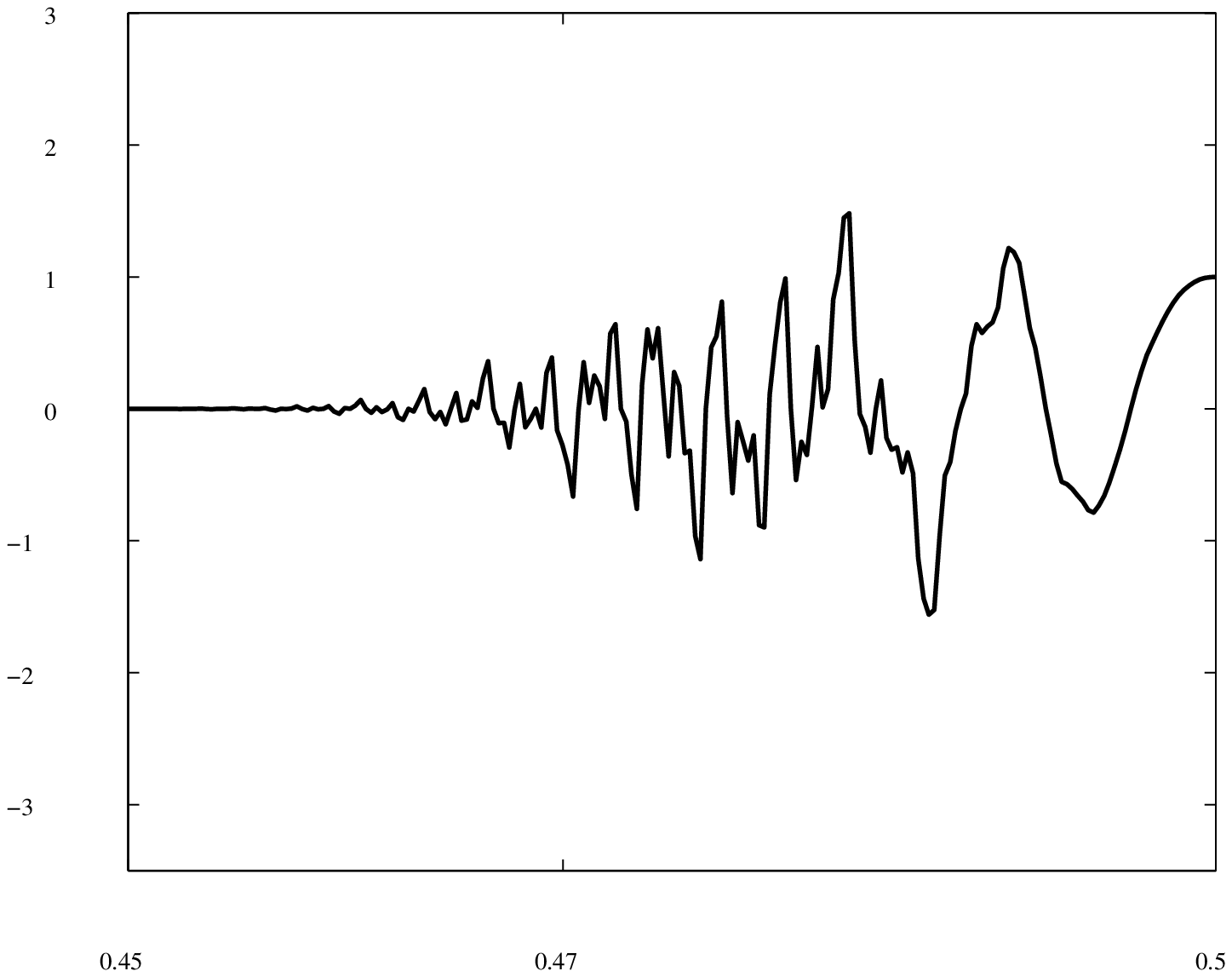}
}
\subfigure[10 point scheme  (zooming)]
{
\psfrag{0.45}[l][l]{\footnotesize{$0.45$ }}
 \psfrag{0.5}[l][l]{\footnotesize{$0.5$ }}
 \psfrag{0.47}[l][l]{\footnotesize{$0.47$ }}
 \psfrag{-3}[l][l]{\footnotesize{$-3$ }}
 \psfrag{-2}[l][l]{\footnotesize{$-2$ }}
 \psfrag{-1}[l][l]{\footnotesize{$-1$ }}
\psfrag{0}[l][l]{\footnotesize{$0$ }}
\psfrag{1}[l][l]{\footnotesize{$1$ }}
\psfrag{2}[l][l]{\footnotesize{$2$ }}
\psfrag{3}[l][l]{\footnotesize{$3$ }}
\includegraphics[width=0.45\textwidth]{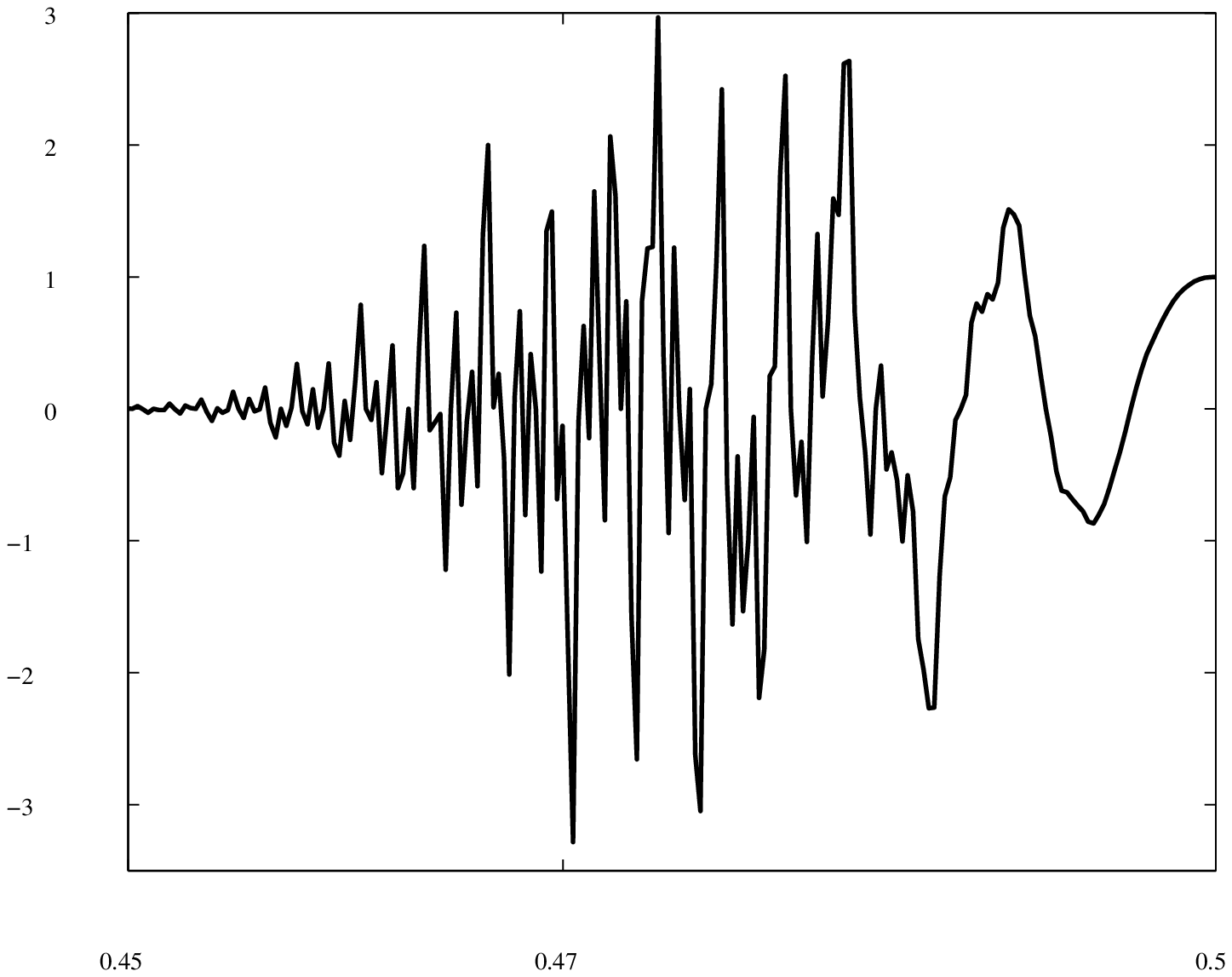}
 }
\caption{Iterated function at scale $13$ starting from $f^8$ defined by $f^8_k=1$ if $k2^{-8}\geq .5$ and $f^8_k=0$ elsewhere}
\label{decentre}
\end{figure}

\indent  In this subsection, we focus on completely decentred Lagrange interpolatory
linear schemes.
In order to apply our theoretical results, we  consider $S$ the two point centered linear scheme and  express any  right hand side excentred scheme $S_{P}$ (where $P$ stands for the number of point of the considered stencil) as a perturbation of it.
Precisely, if we write $S_{P}(f^j)_{2n+1}=S(f^j)_{2n+1} +F_P(\delta_P f^j)_{2n+1}$  we get:

 If  $P$ is  even,
\[
\begin{array}{llll}
F_P(\delta_P f^j)_{2n+1}=& + & \displaystyle{\sum_{k=2 \ k \, even}^{P-2}}
D^{\frac{k}{2}}f_{n+\frac{k}{2}+1}& \frac{(2k-1)!}{2^{2k}(k-1)!(k+1)!}\\
 & - & \displaystyle{\sum_{k=1
\ k\,
odd}^{P-3}}D^{\frac{k+1}{2}}f_{n+\frac{k+1}{2}} & \frac{(4k+5)(2k-1)!}{2^{2k+1}(k-1)!(k+2)!},
\end{array}
\]
 and,if  $P$  is odd\\
\[
\begin{array}{llll}
F_P(\delta_P f^j)_{2n+1}=& + & \displaystyle{\sum_{k=2 \ k\, even}^{P-3}}
D^{\frac{k}{2}}f_{n+\frac{k}{2}+1} & \frac{(2k-1)!}{2^{2k}(k-1)!(k+1)!}\\
 & - & \displaystyle{\sum_{k=1
\ k\,
odd}^{P-4}}D^{\frac{k+1}{2}}f_{n+\frac{k+1}{2}} & \frac{(4k+5)(2k-1)!}{2^{2k+1}(k-1)!(k+2)!}\\
  & - &
 D^{\frac{P-1}{2}}f_{n+\frac{N-1}{2}} & \frac{(2N-3)!}{2^{2(N-2)}(N-3)!(N-1)!}.
\end{array}
\]
Numerical evaluation of the perturbation terms  for $4\leq P \leq 9$ is given
in table \ref{tabsub}.
\begin{table}[!h]
\begin{center}
\begin{tabular}{c|l}
\hline
 P &  $F_P(Df)_{2n+1}$ \\
\hline
$4$ &\rule[-0.3cm]{0pt}{0.7cm} $\frac{-3}{16}Df_{n+1}+\frac{1}{16}Df_{n+2}$ \\
$5$ &\rule[-0.3cm]{0pt}{0.7cm} $\frac{-3}{16}Df_{n+1}+\frac{1}{16}Df_{n+2}-\frac{5}{128}D^2f_{n+2}$ \\
$6$ & \rule[-0.3cm]{0pt}{0.7cm}$\frac{-3}{16}Df_{n+1}+\frac{1}{16}Df_{n+2}-\frac{17}{256}D^2f_{n+2}+\frac{7}{256}D^2f_{n+3}$ \\
$7$ &\rule[-0.3cm]{0pt}{0.7cm} $\frac{-3}{16}Df_{n+1}+\frac{1}{16}Df_{n+2}-\frac{17}{256}D^2f_{n+2}+\frac{7}{256}D^2f_{n+3}-\frac{21}{1024}D^3f_{n+3}$  \\
$8$ &\rule[-0.3cm]{0pt}{0.7cm}
$\frac{-3}{16}Df_{n+1}+\frac{1}{16}Df_{n+2}-\frac{17}{256}D^2f_{n+2}+\frac{7}{256}D^2f_{n+3}-\frac{75}{2048}D^3f_{n+3}+\frac{33}{2048}D^3f_{n+4}$  \\
$9$ &\rule[-0.3cm]{0pt}{0.7cm}
$\frac{-3}{16}Df_{n+1}+\frac{1}{16}Df_{n+2}-\frac{17}{256}D^2f_{n+2}+\frac{7}{256}D^2f_{n+3}-\frac{75}{2048}D^3f_{n+3}+\frac{33}{2048}D^3f_{n+4}$  \\
 &\rule[-0.3cm]{0pt}{0.7cm} $-\frac{429}{32768}D^4f_{n+4}$\\
\hline
\end{tabular}
\caption{Perturbation term $F_P(Df)_{2n+1}$ for different values of $P$
}
\label{tabsub}
\end{center}
\end{table}
\vspace{1cm}
\ \\
 It also appears that  $S_P$ can be written  naturally as a
 perturbation of $S_{P-2}$, for even values of $P$ and as a perturbation of $S_{P-1}$ for odd values of $P$. Indeed, we have:

 When $P$ is even:

\begin{equation}
\label{pertP-2}
\begin{array}{llllll}
S_P(f)_{2n+1} & =& S_{P-2}(f)_{2n+1}&+&\frac{(2P-3)!}{2^{2(P-2)}(P-3)!(P-1)!} & D^{\frac{P-2}{2}}f_{n+\frac{P}{2}}
\\
& &  &  -  & \frac{(4P-8)(2P-7)!}{2^{2P-5}(P-4)!(P-1)!} & D^{\frac{P-2}{2}}f_{n+\frac{P}{2}-1},
 \end{array}
 \end{equation}
and  when $P$ is odd:
\begin{equation}
\label{pertP-1}
\begin{array}{llllll}
S_P(f)_{2n+1} & = & S_{P-1}(f)_{2n+1} & - &
\frac{(2P-3)!}{2^{2(P-2)}(P-3)!(P-1)!} & D^{\frac{P-1}{2}}f_{n+\frac{P-1}{2}}.
\end{array}
 \end{equation}
In  both cases, its is  easy to check that the function  $F$ defined by $F=F(D^{\frac{P-2}{2}}f)$ when $P$ is  even
 and by  $F=F(D^{\frac{P-1}{2}}f)$  when $P$ is odd, is  linear and continuous.
Therefore,   the convergence can be reached as soon as
the contractivity hypothesis (\ref{h2}) for $
D^{\frac{P-2}{2}}$ or $ D^{\frac{P-1}{2}}$ is satisfied. \\
\ \\
Direct calculations provide the estimates gathered in table \ref{table2}.
\begin{table}[!h]
\begin{center}
\begin{tabular}{c|c|c}
\hline
P& perturbation term& contractivity estimate\\
\hline
4&\rule[-0.3cm]{0pt}{0.7cm}$ F(Df)=-\frac{3}{16}Df_{n+1}+\frac{1}{16}Df_{n+2}$&$
||D(S_4f)||_{\infty}\leq\frac{1}{2}||Df||_{\infty}$\\
5&\rule[-0.3cm]{0pt}{0.7cm}$ F(D^2f)=-\frac{5}{128}D^2f_{n+2}$&$ ||D^2(S_5f)||_{\infty}\leq\frac{1}{2}||D^2f||_{\infty}$\\
6&\rule[-0.3cm]{0pt}{0.7cm}$F(D^2f)=-\frac{17}{256}D^2f_{n+2}+\frac{7}{256}D^2f_{n+3}$&$ ||D^2(S_6f)||_{\infty}\leq
\frac{87}{128}||D^2f||_{\infty}$\\
7&\rule[-0.3cm]{0pt}{0.7cm}$F(D^3f)=-\frac{21}{1024}D^3f_{n+3} $&$ ||D^3(S_7f)||_{\infty}\leq
\frac{367}{512}||D^3f||_{\infty}$\\
8&\rule[-0.3cm]{0pt}{0.7cm}$F(D^3f)=-\frac{75}{2048}D^3f_{n+3}+\frac{33}{2048}D^3f_{n+4} $&$ ||D^3(S_8f)||_{\infty}\leq
\frac{475}{512}||D^3f||_{\infty}.$\\
9&\rule[-0.3cm]{0pt}{0.7cm}$F(D^4f)=-\frac{429}{32768}D^4f_{n+4} $&$ ||D^4(S_9f)||_{\infty}\leq
\frac{54734}{32768}||D^4f||_{\infty}.$\\
\hline
\end{tabular}
\caption{Perturbation term  (see \ref{pertP-2} and \ref{pertP-2}) and contractivity estimate for different values of $P$}
\label{table2}
\end{center}
\end{table}
It then follows from theorem \ref{th1} that all the fully excentred interpolatory
subdivision schemes for   $P \leq 8$ points converge.\\
The following comments can be made for the other situations:

Figure \ref{decentre} represents the completely decentred 9 and 10 points scheme iterated functions at scale $8$ starting from $f^0$. From the zooming in the oscillating region, one can guess that the $9$ point scheme converges while the $10$ point doesn't. In fact, following \cite{Dyn}, and using the so called iterative formalism one observes numerically that the spectral radius of the iterated matrix for $10$ points overshoots the critical value $1$ while the spectral radius of the iterated matrix for $9$ points doesn't, that confirms the guess.

Obviously, for these linear schemes, theorem \ref{th5} applies as soon as theorem \ref{th1} does. Therefore stability is ensured for $P\leq 8$.

\ \\

\ \\
\subsection{ The    6 points WENO subdivision scheme}

\indent  $WENO$ subdivision schemes \cite{LOC}, are constructed using  convex combination of
different interpolatory  polynomials of fixed degree.
For degree $3$ and therefore a 6 point stencil, the $WENO-6$ subdivision is given by:
\begin{eqnarray*}
S_{weno}(f^{j})_{2n+1}&=&\frac{\alpha_2}{16} f^{j}_{n-2}
-\frac{5\alpha_2+\alpha_1}{16} f^{j}_{n-1}
+(1+\frac{5\alpha_2+2\alpha_1}{8} )f^{j}_{n}\\
&&+(1+\frac{5\alpha_0+2\alpha_1}{8} )f^{j}_{n+1}
-\frac{5\alpha_0+\alpha_1}{16} f^{j}_{n+2} +\frac{\alpha_0}{16}
f^{j}_{n+3}
\end{eqnarray*}
where the coefficients $\alpha_{i}$ control the convex combination and therefore
satisfy $\alpha_{i} \geq 0$ and $\alpha_0 +\alpha_1+\alpha_2=1$.\\
\ \\
In \cite{LOC}, these coefficients are defined as:
$$
\alpha_i = \frac{a_i}{ a_0+a_1+a_2}
$$
with
$$
a_i=\frac{d_i}{(\epsilon+b_i)^2}
$$
where $b_i$, defined as a function of the first difference $df$  is an indicator of smoothness while  $d_i$
and $\epsilon $  are fixed positive constants.  A set of
possible values for  these constants is suggested in \cite{JS}.\\
\ \\
 \indent The convergence
of the associated subdivision scheme has been studied in  \cite{CDM}.
We present an alternative proof, using   theorem \ref{th1}.\\
\indent First,  the  $WENO-6$ subdivision scheme is written as a
perturbation of the linear two point interpolation scheme as:
\begin{eqnarray}
\label{weno}
S_{weno}(f^{j})_{2n+1}&=& \frac{f^{j}_n+f^{j}_{n+1}}{2}+\frac{\alpha_0}{16}Df^j_{n+2}\\
\nonumber
         &  &-  \frac{3\alpha_0+\alpha_1}{16}  Df^j_{n+1}-
\frac{\alpha_1+3\alpha_2}{16}  Df^j_{n}
+\frac{\alpha_2}{16}Df^j_{n-1} \\
\nonumber
\end{eqnarray}
with $\alpha_0=\alpha_0(df^j)$, $\alpha_1=\alpha_1(df^j)$ and
$\alpha_2=\alpha_2(df^j)$.\\
\ \\
\ \\
We then have the following proposition\,:\\

\begin{proposition}
The $WENO-6$ subdivision scheme is convergent and, for any initial sequence $f^j$, the
limit function  belongs to  $C^{-log_2(\frac{3}{4})^-}$.
\end{proposition}
\ \\
{\bf Proof}\\

According to remark \ref{deuxoperateurs}, the proof can be performed in three steps considering that $F$ is a function of $df^j$ and $Df^j$.

First, according to the definition of $F$ and to the properties of  $\alpha_i$, we have
\begin{eqnarray*}
 |F(d,D)| &\leq & \frac{1}{2}\max{\left  ( ||d||_{\infty},||D||)_{\infty}\right ). }
 \end{eqnarray*}

Second, we  prove (\ref{deux1}) for the first difference operator  $d$:

We have, for $f \in l^{\infty}$\,:
\[ d(S_{weno}(f))_{k}=S_{weno}(f)_{k+1}-S_{weno}(f)_{k}\]
\ \\
We have to consider two cases, according to the parity of $k$. We give the details for $k=2n+1$, the even case being similar.

\begin{eqnarray*}
d(S_{weno}(f))_{2n+1} & = & S_{weno}(f)_{2n+2}-S_{weno}(f)_{2n+1}\\
 & = & f_{n+1}-\frac{f_n+f_{n+1}}{2} -\frac{\alpha_0}{16}Df_{n+2}+  \frac{3\alpha_0+\alpha_1}{16}  Df_{n+1}\\
& & +
\frac{\alpha_1+3\alpha_2}{16}  Df_{n}
-\frac{\alpha_2}{16}Df_{n-1} \\
 &=& \frac{ df_n}{2}  -\frac{\alpha_0}{16}Df_{n+2}+
\frac{3\alpha_0+\alpha_1}{16} Df_{n+1}\\
& &+
\frac{\alpha_1+3\alpha_2}{16}  Df_{n}
-\frac{\alpha_2}{16}Df_{n-1} \\
\end{eqnarray*}
Since   $\alpha_0+\alpha_1+\alpha_2=1 $ and $0<\alpha_1<1$, we have:
\begin{eqnarray}
\nonumber
 |d(S_{weno}(f))_{2n+1}|& \leq &\frac{1}{2}||df||_{\infty}+\left ( \frac{\alpha_0}{16}+ \frac{3\alpha_0}{16}
+\frac{\alpha_1}{16}+
\frac{\alpha_1}{16}+\frac{3\alpha_2}{16} +\frac{\alpha_2}{16}\right ) || Df||_{\infty}, \\
\nonumber
& \leq &\frac{1}{2}||df||_{\infty}+\frac{4-2\alpha_1}{16} || Df||_{\infty},
\\
\nonumber
& \leq & \frac{1}{2}||df||_{\infty}+\frac{1}{4} || Df||_{\infty},
\\
 \label{eq23}
& \leq &  \frac{3}{4} \max{\left ( || df||_{\infty},||Df||_{\infty}. \right ) }
\end{eqnarray}

Third, we  prove  inequality (\ref{h2}) for the second difference operator $D$.\\
Again, two cases have to be  considered:
\ \\
$\bullet$  For k=2n+1, then\,:\\
\begin{eqnarray*}
D(S_{weno}(f))_{2n+1}&=&f_{n+1}-2S_{weno}(f)_{2n+1} +f_n,\\
&=&  -\frac{\alpha_0}{8}Df_{n+2}+\frac{3\alpha_0+\alpha_1}{8} Df_{n+1}+
\frac{\alpha_1+3\alpha_2}{8} Df_{n}-\frac{\alpha_2}{8}Df_{n-1}.\\
\end{eqnarray*}
Using $0<\alpha_0<1 $, $0<\alpha_1<1$, $0<\alpha_2<1 $  and $\alpha_0+\alpha_1+\alpha_2=1$, we get:
\begin{eqnarray}
  \nonumber
 |D(S_{weno}(f))_{2n+1}|& \leq &  \frac{4\alpha_0+2\alpha_1+4\alpha_2}{8}  \Vert Df \Vert_{\infty},\\
\nonumber
 & \leq& \frac{4-2\alpha_1}{8}  \Vert Df \Vert_{\infty},\\
\label{eq24}
  & \leq & \frac{1}{2} \Vert Df \Vert_{\infty}.
      \end{eqnarray}
\ \\
\ \\
$\bullet$ For k=2n, then\,:\\
\begin{eqnarray*}
D(S_{weno}(f))_{2n}&=& S_{weno}(f)_{2n+1}-2f_n+S_{weno}(f)_{2n-1},\\
 &=&  \frac{f_{n+1}-2f_n+f_{n-1}}{2}+  \frac{\alpha_0}{16} Df_{n+2}\\
  & & - \frac{2\alpha_0+\alpha_1}{16} Df_{n+1}
  -  \frac{3\alpha_0+2\alpha_1}{16} + \frac{3\alpha_2}{16} Df_{n}\\
  & &-\frac{\alpha_1+2\alpha_2}{16}Df_{n-1}+\frac{\alpha_2}{16}Df_{n-2}.\\
\end{eqnarray*}
Using $\alpha_0+\alpha_1+\alpha_2=1$, we get:\\
\begin{eqnarray*}
D(S_{weno}(f))_{2n}& =& \frac{5}{16} Df_{n}
 + \frac{\alpha_0}{16}Df_{n+2}-
 \frac{2\alpha_0+\alpha_1}{16}
 Df_{n+1}+\frac{\alpha_1}{16}Df_{n}\\
  & & -  \frac{\alpha_1+2\alpha_2}{16}  Df_{n-1}+\frac{\alpha_2}{16}Df_{n-2}\\
\end{eqnarray*}
Then, with $0<\alpha_0<1 $, $0<\alpha_1<1$ and  $0<\alpha_2<1 $\,:
\begin{eqnarray}
\nonumber
 |D(S_{weno}(f))_{2n}|& \leq & \left ( \frac{5}{16}+\frac{3\alpha_0+3\alpha_1+3\alpha_2}{16}\right )  \Vert Df \Vert_{\infty},\\
 \nonumber
& \leq & \left ( \frac{5}{16}+\frac{3}{16}\right ) \Vert Df \Vert_{\infty},\\
\label{eq25}
& \leq & \frac{1}{2} \Vert Df \Vert_{\infty}.
\end{eqnarray}
\ \\
Therefore, from (\ref{eq23}), (\ref{eq24}) and (\ref{eq25}), we obtain the
inequality\,:
\begin{eqnarray*}
\exists c<1 \quad \forall f \in l^{\infty} \qquad
\max{\left (||d(S_{weno}(f))||_{\infty},||D(S_{weno}(f))||_{\infty} \right
)}\leq \frac{3}{4} \max{\left ( || df||_{\infty},||Df||_{\infty} \right ) }.
\end{eqnarray*}

Finally, using  remarks \ref{deuxoperateurs} and \ref{regularite}
we get  the convergence of the $WENO-6$ subdivision scheme  to a
$C^{-log_2(\frac{3}{4})^-}$ function.

$\Box$\\
\ \\
\subsection{Power-P subdivision scheme: definition and convergence}
\ \\
\indent In the same vein as  the PPH scheme (\cite{ADLT}), the power P scheme is a four point scheme
based on a
piecewise degree $3$  polynomial prediction. Considering $S_{\cal{L}}$ the centered four point Lagrange interpolation prediction that reads:

\begin{equation}\label{predic3}
(S_{\cal{L}}(f))_{2n+1}= \frac{f_n+f_{n+1}}{2}-\frac{1}{8}
\frac{ Df_{n+1}+ Df_{n}}{2},
\end{equation}

the definition of the Power-P subdivision scheme is based on the substitution of the  arithmetic mean of second
order differences, $\frac{ Df_{j+1}+ Df_{j}}{2}$,
by a general mean $power_p( Df_{j},Df_{j+1})$ defined in \cite{SM} for any integer $p\geq 1$, and any couple $(x,y)$  as:

\begin{eqnarray}
\label{defpower}
power_p(x,y) & = & \frac{sign(x)+sign(y)}{2}\frac{x+y}{2}\left ( 1- \left |
\frac{x-y}{x+y} \right | ^p \right ).
\end{eqnarray}
Note that it coincides for  $p=1$, with  the arithmetic mean and for $p=2$ with  the geometric mean.The Power-P subdivision scheme then  naturally appears as a perturbation of the linear two point interpolation scheme since it is defined by
\begin{eqnarray}
S_{power_p}(f^j])_{2n+1}&=&\frac{f_n+f_{n+1}}{2}-\displaystyle{\frac{1}{8}}power_p(Df^j_n,Df^j_{n+1}).
\end{eqnarray}

Before establishing the convergence theorem we first prove the following lemma:

\begin{lemma}
\label{lem1}
For any $(x,y),(x',y') \in \mathbb{R}^2$, the function $power_p $ satisfies the
following properties\,:
\begin{enumerate}
\item $power_p(x,y)=power_p(y,x)$
\item $power_p(x,y)=0 \quad if\; xy\leq 0$
\item $power_p(-x,-y)=-power_p(x,y)$
\item $|power_p(x,y)|\leq \max{(|x|,|y|)}$
\item $|power_p(x,y)| \leq p \min{(|x|,|y|)}$
\end{enumerate}
\end{lemma}
\ \\
{\bf Proof}\\
\ \\
Claims of $1-4$ are obvious;\\
Inequality  $5$ comes from the equality
\begin{eqnarray*}
power_p(x,y) & = & \frac{sign(x)+sign(y)}{2} min(x,y)\left [
1+\left | \frac{x-y}{x+y} \right | +\dots  +\left | \frac{x-y}{x+y}\right |^{p-1}\right ].
\end{eqnarray*}
$\Box$\\
 We then have the following proposition:

\begin{proposition}
The Power P subdivision scheme is uniformly convergent and, for any initial sequence $f^j$ the limit function belongs to
$C^{1^-}$ for $p\leq 4$ and $C^{-log_2(\frac{3}{4})^-}$ for $p\geq 5 $.
\end{proposition}
\ \\
{\bf Proof}\\

Here again,  the hypotheses of the general theorem \ref{th1} must be checked:

We first check hypothesis (\ref{h1}).  Using property $4$ of  lemma \ref{lem1}, we obtain for
$d\in l^{\infty}$\,:
\begin{eqnarray*}
 |F(d)| &\leq & \frac{1}{8}\max{\left(|d_n|,|d_{n+1}|\right ) }\\
 |F(d)| &\leq & \frac{1}{8} ||d ||_{\infty}\\
 \end{eqnarray*}
 \ \\

 Then we consider hypothesis (\ref{h2}):

We study as before  two different cases\,:\\
\ \\
$\bullet$ For $k=2n+1$:
\begin{eqnarray*}
D(S_{power_p}(f))_{2n+1}&=&=f_n -2S_{power p}(f)_{2n+1}+ f_{n+1}\\
&=& f_{n+1}+f_{n}-2\frac{f_n+f_{n+1}}{2}
+2\frac{1}{8}power_p(Df_n,Df_{n+1})\\
&=&  \frac{1}{4}power_p(Df_n,Df_{n+1})
\end{eqnarray*}
\ \\
 From  property $4$ of  lemma \ref{lem1} we get:
\begin{eqnarray}
 \label{power1}
|D(S_{power_p}(f))_{2n+1}| & \leq & \frac{1}{4} \Vert Df\Vert_{\infty}.
\end{eqnarray}
\ \\
\ \\
$\bullet$ For $k=2n$:
\begin{eqnarray*}
D(S_{power_p}(f))_{2n} & = & S_{power_p}(f)_{2n-1}-2f_{n}+S_{power
p}(f)_{2n+1} \\
&=&
\frac{f_n+f_{n+1}}{2}-\frac{1}{8} power_p (Df_n,Df_{n+1})-2f_{n}\\
& & +\frac{f_{n-1}^j+f_{n}}{2}-\frac{1}{8}power_p(Df_{n-1},Df_{n})
\\
& = & \frac{Df_n}{2}-\frac{1}{8} \left
 ( power_p(Df_n,Df_{n+1})+power_p(Df_{n-1},Df_{n}) \right ) \\
\end{eqnarray*}
\ \\
For $p\geq 5$, from property $4$ of  lemma \ref{lem1} we get:
\begin{eqnarray}
\label{power2}
|D(S_{power_p}(f))_{2n}|  &\leq & \frac{3}{4} ||Df||_{\infty}.
\end{eqnarray}
\ \\
\ \\
For $p\leq 4$, we note $D(S_{power_p}(f))_{2n}=Z(Df_n,Df_{n+1},Df_{n-1})$ with
\[ Z(x,y,z)=\frac{x}{2}-\frac{1}{8}(power_p(x,y)+power_p(x,z))\]
From definition \ref{defpower} and property $4$ and $5$ of lemma \ref{lem1}, we have,\\
\ \\
if $x>0$,
\begin{eqnarray*}
\frac{x}{2}-\frac{1}{8}(\max{(x,y)}+\max{(x,z)})  \leq & Z(x,y,z) & \leq \frac{x}{2}\\
\frac{x}{4} \leq & Z(x,y,z) & \leq \frac{x}{2}\\
0 \leq  & |Z(x,y,z)| & \leq \frac{1}{2}|x|\\
\end{eqnarray*}
if $x< 0$,
\begin{eqnarray*}
 \frac{x}{2} \leq & Z(x,y,z) & \leq
 \frac{x}{2}+\frac{p}{8}(\min{(|x|,|y|)}+\min{(|x|,|z|)}) \\
\frac{x}{2} \leq  & Z(x,y,z) & \leq (\frac{p}{4}-\frac{1}{2})|x|\\
0 \leq & |Z(x,y,z)| & \leq \frac{1}{2}|x|\\
\end{eqnarray*}
Finally,
\begin{equation}
\label{power2'}
| D(S_{power_p}(f))_{2n}|  \leq \frac{1}{2} \Vert Df \Vert_\infty \qquad
.
\end{equation}
\ \\
From (\ref{power1}), (\ref{power2}) and (\ref{power2'}), we obtain\,:
\begin{eqnarray*}
\Vert DS_{power_p}(f)\Vert_\infty & \leq &\frac{1}{2} \Vert Df \Vert_\infty \qquad
for \; p \leq 4.\\
\Vert DS_{power_p}(f)\Vert_\infty & \leq &\frac{3}{4} \Vert Df \Vert_\infty \qquad
for \; p \geq 5.\\
\end{eqnarray*}
Finally,  theorem \ref{th1} and remark \ref{regularite}
provides  the convergence to a
$C^{1^-}$ if $p\leq 4$ and $C^{-log_2(\frac{3}{4})^-}(\mathbb{R})$ if $p\geq 5$.\\
$\Box$
\ \\
\subsection{The convergence of a non linear scheme using spherical coordinates}

\indent The non linear subdivision scheme studied in this section  is defined in
\cite{AEV} where it is considered as a non regular  interpolatory  subdivision scheme using local spherical coordinates. Here, we consider it as a regular  subdivision scheme applied to the $I\!\!R^2$ point  sequence  $P_n^j(x_n^j, f_n^j)^t_{n \in Z\!\!Z}$.
The resulting scheme reads (see \cite{AEV}):

\begin{eqnarray}
\label{spherical}
\left (\begin{array}{cc}
x_{2n+1}^{j+1}\\
f_{2n+1}^{j+1}
\end{array}\right )=\left (\begin{array}{cc}
\frac{x^j_{n}+x^j_{n+1}}{2}\\
\frac{f^j_{n}+f^j_{n+1}}{2}
\end{array}\right )
+\frac{r_n^j}{4}\left  ( \begin{array}{cc}
 cos \left ( \theta_n^j+h(\alpha_n^j) \right ) - cos  \left (\theta_{n+1}^j+h(\beta_{n+1}^j)\right )  \\
  sin \left ( \theta_n^j+h(\alpha_n^j) \right ) - sin \left (\theta_{n+1}^j+h(\beta_{n+1}^j)\right )
  \end{array} \right )
  \end{eqnarray}

with\,:
\begin{eqnarray}
\label{r}
r^j_{n}& =& \sqrt{(x^j_{n+1}-x^j_{n})^2+(f^j_{n+1}-f^j_{n})^2}, \\
\label{theta}
\theta^j_{n}& =& \arctan \left (\frac{f^j_{n+1}-f^j_{n-1}}{x^j_{n+1}-x^j_{n-1}} \right ), \\
\label{gamma}
\gamma^j_n & =& \arctan \left (\frac{f^j_{n+1}-f^j_{n}}{x^j_{n+1}-x^j_{n}} \right ), \\
\label{alpha}
\alpha^j_n & = &\gamma_n^j-\theta^j_n,\\
\label{beta}
\beta_{n+1}^j &=&\gamma_{n}^{j}-\theta^j_{n+1},
\end{eqnarray}
and, $\theta^j_n$, $\gamma^j_n$ $\in [-\frac{\pi}{2};-\frac{\pi}{2} ]$.\\
As explained in \cite{AEV}, the design  of $h$, is performed to produce regular limit  functions. It is then defined as  as a $C^1$ function that is contractive for small values of $\alpha$ and that coincides with identity for large value of $\alpha$. Note that $h=0$ provides   the classical linear two point centered scheme.

\ \\
In our context, we will note this  scheme $S_{spherical}$, and $S_1$, $S_2$
will stand for the  schemes associated to each  coordinates\,:
We then get
\begin{eqnarray*}
S_1(x,f) _{2n+1}& = & \frac{x_{n}+x_{n+1}}{2}+(F_1(dx,df))_{2n+1}, \\
S_2(x,f)_{2n+1} & = & \frac{f_{n}+f_{n+1}}{2}+(F_2(dx,df))_{2n+1}, \\
\end{eqnarray*}
with
\begin{eqnarray*}
(F_1(dx,df))_{2n+1}& =& \frac{r_n^j}{4}\left ( cos \left ( \theta_n^j+h(\alpha_n^j)
\right ) - cos  \left (\theta_{n+1}^j+h(\beta_{n+1}^j)\right ) \right ), \\
(F_2(dx,df))_{2n+1}&=& \frac{r_n^j}{4}\left  (sin \left ( \theta_n^j+h(\alpha_n^j)
\right ) -sin \left (\theta_{n+1}^j+h(\beta_{n+1}^j)\right ) \right ).  \\
\end{eqnarray*}
From (\ref{r}, \ref{theta}, \ref{gamma}),
$r_n$, $\theta_n$ and $\gamma_n$ can be
written using the first divided difference  $(df^j,dx^j)$ as:
\begin{eqnarray*}
r_n & = & \sqrt{(dx^j_n)^2+(df^j_n)^2}\\
\theta_n & = & \arctan \left( \frac{df^j_{n}+df^j_{n-1}}{dx^j_{n}+dx^j_{n-1}} \right )\\
\gamma_n & =& \arctan \left ( \frac{df^j_n}{dx^j_n}\right )
\end{eqnarray*}
as well as $\alpha_n$ and $\beta_n$ thanks to  (\ref{alpha}) and
(\ref{beta}).\\
We then have the following proposition:
\begin{proposition}
The scheme $S_{spherical}$ defined in  (\ref{spherical})  is convergent.
\end{proposition}
{\bf Proof}\\
\ \\
We again check the hypotheses of theorem \ref{th1} generalized to $I\!\!R^2$ according to remark  \ref{R2}.
 We have,
\begin{equation}
\label{prop r_n}
r_n \leq \sqrt{2} \max{\left ( |dx_n|, |df_n|\right )},
\end{equation}
 and therefore, for $i=1,2$:
\begin{eqnarray*}
|(F_i(dx,df))_{2n+1}| & \leq  & 2 \frac{\sqrt{2}\max{\left ( |dx_n|, |df_n| \right ) }}{4},\\
  & \leq  & \frac{\sqrt{2}}{2} \max{\left ( ||dx||_{\infty},||df||_{\infty}\right )},\\
\end{eqnarray*}
that shows  that  the hypothesis  (\ref{h1}) of theorem \ref{th1} is satisfied.\\
\ \\
\ \\
We now check  hypothesis (\ref{h2}).
\ \\
For  $f \in l^{\infty}$ we have
\begin{eqnarray*}
d(S_1(x,f))_{2n}& = & S_1(x,f)_{2n+1}-S_1(x,f)_{2n}\\
 & = &  \frac{x_{n}+x_{n+1}}{2}+\frac{r_n}{4}\left ( cos \left ( \theta_n+h(\alpha_n) \right ) - cos  \left (\theta_{n+1}+h(\beta_{n+1})\right )  \right )- x_n\\
 & = &  \frac{x_{n+1}-x_n}{2}+\frac{r_n}{4}\left ( cos \left ( \theta_n+h(\alpha_n) \right ) - cos  \left (\theta_{n+1}+h(\beta_{n+1})\right )  \right ),\\
\end{eqnarray*}
and therefore
\begin{eqnarray*}
|d(S_1(x,f))_{2n}| &\leq  & \frac{|| dx||_{\infty}}{2}+\frac{\sqrt{2}\max{\left ( ||dx||_{\infty}, ||df||_{\infty} \right )}}{4}\left |  \theta_n+h(\alpha_n)  -\theta_{n+1}-h(\beta_{n+1}) \right |.\\
\end{eqnarray*}
Using the definitions of  $\alpha_n$ and  $\beta_n$ we get

\begin{eqnarray*}
|d(S_1(x,f))_{2n}|& \leq &
\frac{|| dx||_{\infty}}{2}+\frac{\sqrt{2}\max{\left ( ||dx||_{\infty}, ||df||_{\infty} \right )}}{4}
\left |\theta_n+h(\gamma_n-\theta_n)-\left
(\theta_{n+1}+h(\gamma_n-\theta_{n+1})\right ) \right |\\
\end{eqnarray*}
 and
\begin{eqnarray*}
|d(S_1(x,f))_{2n}| &  \leq&  \frac{ || dx||_{\infty}}{2}+\frac{\sqrt{2}\max{\left (
||dx||_{\infty}, ||df||_{\infty} \right )}}{4} \max_{x\in [-\pi, \pi]}{\left
(1-h'(x) \right )} \left |  \theta_n-\theta_{n+1} \right |\\
 &  \leq &  \left ( \frac{1}{2}+\frac{\sqrt{2}\pi}{4} \max_{x\in [-\pi,
 \pi]}{|1-h'(x)|} \right ) \max{\left ( ||dx||_{\infty}, ||df||_{\infty} \right )}\\
\end{eqnarray*}

The contractivity hypothesis  (\ref{h2}) is therefore satisfied as soon as  $ \forall x \in [-\pi,\pi], 1-\frac{\sqrt{2}}{\pi}< h'(x)<1+\frac{\sqrt{2}}{\pi}$. For instance, the following function h:

\begin{equation*}
h(x)= \left \{
\begin{array}{lllllll}
x & if & -\pi & < & x  & \leq  & -\frac{\pi}{2} \\
x+\frac{63}{125 \pi}(x+\frac{\pi}{2})2- \frac{3969}{625 \pi ^2}(x+\frac{\pi}{2})(x+\frac{\pi}{7} )  & if & -\frac{\pi}{2} & < &x &< &-\frac{\pi}{7}\\
0.55x & if  & -\frac{\pi}{7} & \leq  & x &  \leq  & \frac{\pi}{7}\\
x-\frac{63}{125 \pi}(x-\frac{\pi}{2})2- \frac{3969}{625 \pi ^2}(x-\frac{\pi}{2})(x-\frac{\pi}{7} )& if  &  \frac{\pi}{7} & < & x & < & \frac{\pi}{2}\\
x & if &  \frac{\pi}{2} & \leq  & x &  <  & \pi \end{array}
\right.
\end{equation*}
which is in agreement with the criteria proposed in  \cite{AEV} leads to a scheme satisfying (\ref{h2}).

Since the same sketch of proof also provides the contractivity for $|d((S_2(x,f))_{2n}|$ we get the convergence applying theorem  \ref{th1}.

$\Box$\\
\section{Conclusion}

We have formulated  convergence and stability conditions for non linear subdivision schemes and associated multi-resolutions. These conditions deal with the difference with a suitable linear and convergent subdivision scheme. Many examples show that this formulation lead to simple proofs of convergence and stability.


\begin{thebibliography}{99}



\bibitem{ADLT} S.Amat, R.Donat, J.Liandrat and J.C.Trillo,
{\em Analysis of a new nonlinear subdivision scheme. Applications
in image processing}, to appear in Foundations of Computational
Mathematics.



\bibitem{AL} S.Amat and J.Liandrat. {\em On the stability
of PPH nonlinear multi-resolution}. Applied and Computational
Harmonic Analysis. {\bf 18} (2), 198-206, (2005).


\bibitem{AD00}  F.Ar\`andiga and R.Donat.
{\em Nonlinear Multi-scale Decompositions: The approach
of A.Harten.}
 Numerical Algorithms, {\bf 23},  175-216, (2000).

\bibitem{B}  G. Beylkin. {\em Wavelets, Multi-resolution Analysis and Fast
Numerical Algorithms.} { INRIA lectures, manuscript, (1991)}.

\bibitem{CDM} A.Cohen, N.Dyn and B.Matei. {\em Quasi linear
subdivision schemes with applications to ENO interpolation}. {
Applied and Computational Harmonic Analysis}, {\bf 15}, 89-116,
(2003).

\bibitem{DeDu} G.Delauriers and S.Dubuc. {\em Symmetric  Iterative Interpolation Processes}.
{ Const. Approx.}, {\bf 5},
49-68, (1989).



\bibitem{DRS} I. Daubechies, O. Runborg  and W. Sweldens.
{\em Normal multi-resolution approximation of curves}, { Const.
Approx.}, 20, pp 399-463, 2004.




\bibitem{DY} D.Donoho and T.P. Yu. {\em Nonlinear pyramid transforms
based on median interpolation}. { SIAM J. Math. Anal.}, {\bf
31}(5),  1030-1061, (2000).


\bibitem{FM} M.S. Floater  and C.A. Michelli. {\em Nonlinear
 stationary subdivision}, { Approximation theory: in memory of A.K. Varna},
  {\em edt: Govil N.K, Mohapatra N., Nashed Z., Sharma A., Szabados J.}, 209-224, (1998).



\bibitem{Har93}  A. Harten.
{\em Discrete Multi-resolution Analysis and Generalized Wavelets.}
{ J. Appl. Numer. Math., {\bf 12}, 153-192, (1993).}



\bibitem{Har96}  A. Harten. {\em Multi-resolution Representation of Data II:
General Framework}. { SIAM J. Numer. Anal. {\bf 33} 3, 1205-1256,
(1996).}

\bibitem{KD} F. Kuijt and R. van Damme. {\em Convexity preserving
interpolatory subdivision schemes}. { Const. Approx.}, {\bf 14},
609-630, (1998).



\bibitem{LOC} X.D.Liu, S.Osher and T.Chan. {\em Weighted
essentially non-oscillatory schemes.} Journal of Computational
Physics, {\bf 115},  200-212, (1994).


\bibitem{Os} P.Oswald. {\em
    Smoothness of Nonlinear Median-Interpolation Subdivision},
      { Adv. Comput. Math.}, {\bf 20}(4), 401-423, (2004).


\bibitem{Dyn} N.Dyn. {\em
    Subdivision schemes in computer aided geometric design},
      { Oxford University Press}, {\bf 20}(4), 36-104, (1992).

\bibitem{SM} S.Serna and A.Marquina. {\em
   power ENO methods: a fifth order accurate Weighted Power ENO method},
      { Journal of Computational Physics}, {\bf 194}, 632-658, (2004).

\bibitem{JS} Guang-Shan Jiang  and Chi-Wang Shu. {\em
  Efficient Implementation of Weighted ENO Schemes},
      { Journal of Computational Physics}, {\bf 126}, 202-228, (1996).

\bibitem{AEV} N.Aspert, T.Ebrahimi  and P.Vandergheynst. {\em
 Non-linear subdivision using local coordinates},
       {Computer Aided Geometric Design}, {\bf 20}, 165-187, (2003).

\end{thebibliography}
\end{document}